\documentclass{amsart}

\usepackage{fancyhdr}
\usepackage{lastpage}
\usepackage{stmaryrd,yhmath}

\newcommand{\bdis}{\begin{displaymath}}
\newcommand{\edis}{\end{displaymath}}
\newcommand{\be}{\begin{equation}}
\newcommand{\ee}{\end{equation}}
\newcommand{\mbb}{\mathbb}
\newcommand{\mcal}{\mathcal}

\newcommand{\vp}{\varphi}

\newcommand{\zf}{\zeta\left(\frac{1}{2}+it\right)}

\theoremstyle{definition}

\newtheorem{cor}[]{Corollary}

\theoremstyle{remark}
\newtheorem{remark}[]{Remark}

\newtheorem*{mydef1}{{\bf Theorem}}

\newtheorem*{mydef5}{{\bf Lemma}}

\numberwithin{equation}{section}

\begin{document}

\title{Jacob's ladders and the nonlocal interaction of the function $Z(t)$ with the function $\tilde{Z}^2(t)$ on the distance
$\sim (1-c)\pi(t)$ for a collection of disconnected sets}

\author{Jan Moser}

\address{Department of Mathematical Analysis and Numerical Mathematics, Comenius University, Mlynska Dolina M105, 842 48 Bratislava, SLOVAKIA}

\email{jan.mozer@fmph.uniba.sk}


\begin{abstract}
It is shown in this paper that there is a fine correlation of the third order between the values of the functions $Z[\vp_1(t)]$ and $\tilde{Z}^2(t)$
which corresponds to two collections of disconnected sets. The corresponding new asymptotic formula cannot be obtained within known theories of
Balasubramanian, Heath-Brown and Ivic.
\end{abstract}

\maketitle

\section{Result}

\subsection{}

Let (see \cite{2}, (3), (4))
\be \label{1.1}
\begin{split}
G_1(x) & = G_1(x;T,H)=\bigcup_{T\leq t_{2\nu}\leq T+H}\left\{ t:\ t_{2\nu}(-x)<t<t_{2\nu}(x),\ 0<x\leq\frac{\pi}{2}\right\} \\
G_2(y) & = G_2(y;T,H)=\\
& = \bigcup_{T\leq t_{2\nu+1}\leq T+H}\left\{ t:\ t_{2\nu+1}(-y)<t<t_{2\nu+1}(y),\ 0<y\leq\frac{\pi}{2}\right\}
\end{split}
\ee
\be \label{1.2}
H=T^{1/6+2\epsilon} ,
\ee
and the collection of the sequences $\{t_\nu(\tau)\},\ \tau\in [-\pi,\pi],\ \nu=1,2,\dots $ be defined by the equation (see \cite{2}, (1))
\be \label{1.3}
\vartheta[t_\nu(\tau)]=\pi\nu+\tau;\ t_\nu(0)=t_\nu .
\ee

\subsection{}

In this paper we obtain some new properties of the signal
\bdis
Z(t)=e^{i\vartheta(t)}\zf
\edis
generated by the Riemann zeta-function. Namely, let
\be \label{1.4}
G_1(x)=\vp_1[\mathring{G}_1(x)],\ G_2(y)=\vp_1[\mathring{G}_2(y)]
\ee
where $y=\vp_1(T),\ T\geq T_0[\vp_1]$ is the Jacob's ladder. The following theorem holds true.

\begin{mydef1}
\be \label{1.5}
\begin{split}
\int_{\mathring{G}_1(x)}Z[\vp_1(t)]\tilde{Z}^2(t){\rm d}t & = \frac 2\pi H\sin x + \mcal{O}(T^{1/6+\epsilon}) , \\
\int_{\mathring{G}_2(y)}Z[\vp_1(t)]\tilde{Z}^2(t){\rm d}t & = -\frac 2\pi H\sin y + \mcal{O}(T^{1/6+\epsilon}) ,
\end{split}
\ee
where
\be \label{1.6}
t-\vp_1(t)\sim (1-c)\pi(t),\qquad t\to \infty ,
\ee
and $c$ is the Euler's constant and $\pi(t)$ is the prime-counting function.
\end{mydef1}

Let us remind another representation of the signal $Z(t)$ given by the Riemann-Siegel formula (in its local form)
\bdis
Z(t)=2\sum_{n<P_0}\frac {1}{\sqrt{n}}\cos\{\vartheta(t)-t\ln n\}+\mcal{O}(T^{-1/4}),\ t\in [T,T+H],\ P_0=\sqrt{\frac{T}{2\pi}} ,
\edis
i.e. in the form of the resultant oscillation of the system of nonlinear oscillators
\bdis
\frac {2}{\sqrt{n}}\cos\{\vartheta(t)-t\ln n\};\ \vartheta(t)=\frac t2\ln\frac {t}{2\pi}-\frac t2-\frac{\pi}{8}+\mcal{O}\left(\frac 1t\right)
\edis
(at the same time see, for example, the expression $Z[t_\nu(\tau)]$ by (\ref{1.3})).

\subsection{}

Let (comp. (\ref{1.4}))
\be \label{1.7}
T=\vp_1(\mathring{T}),\ T+H=\vp_1(\widering{T+H}) .
\ee
Since (see (\ref{1.6}), (\ref{1.7}))
\bdis
\mathring{T}-\vp_1(\mathring{T})\sim (1-c)\frac{\mathring{T}}{\ln\mathring{T}} \ \Rightarrow \ \mathring{T}-T\sim (1-c)\frac{\mathring{T}}
{\ln\mathring{T}} \ \Rightarrow \ \mathring{T}\sim T
\edis
we have seen (see (\ref{1.2}))
\bdis
\mathring{T}-(T+H)\sim (1-c)\frac{\mathring{T}}{\ln\mathring{T}}-H\sim (1-c)\frac{\mathring{T}}{\ln\mathring{T}} ,
\edis
i.e. $\mathring{T}>T+H$. Then we have
\be \label{1.8}
\begin{split}
& [T,T+H]\cap [\mathring{T},\widering{T+H}]=\emptyset;\ T+H<\mathring{T} , \\
& \rho\{[T,T+H];[\mathring{T},\widering{T+H}]\}\sim (1-c)\pi(T) ,
\end{split}
\ee
where $\rho$ stands for the distance of the corresponding segments (comp. \cite{5}, (1.3), (1.6)).

\begin{remark}
Some nonlocal interaction of the functions $Z[\vp_1(t)]$ and $\tilde{Z}^2(t)$ is expressed by the formula (\ref{1.5}) where (see (\ref{1.7}))
\bdis
t\in\mathring{G}_1(x)\cup \mathring{G}_2(y)\cap [\mathring{T},\widering{T+H}] \ \Rightarrow \ \vp_1(t)\in G_1(x)\cup G_2(y)\cap [T,T+H] .
\edis
Such an interaction is connected with two collections of disconnected sets unboundedly receding each from other (see (\ref{1.8}), $\rho\to\infty$ as
$T\to\infty$) - like mutually receding galaxies (the Hubble law). Compare this remark with the Remark 3 in \cite{5}.
\end{remark}

\begin{remark}
The asymptotic formulae (\ref{1.5}) cannot be obtained by methods of Balasubramanian, Heath-Brown and Ivic (comp. \cite{1}).
\end{remark}

This paper is a continuation of the series \cite{3}-\cite{15}.

\section{The first corollaries}

First of all, we obtain from (\ref{1.5})

\begin{cor}
\be \label{2.1}
\int_{\mathring{G}_1(x)\cup\mathring{G}_2(y)}Z[\vp_1(t)]\tilde{Z}^2(t){\rm d}t=\left\{\begin{array}{lcr}
\frac{2}{\pi}(\sin x-\sin y)H+\mcal{O}(T^{1/6+\epsilon}) & , & x\not= y \\
\mcal{O}(T^{1/6+\epsilon}) & , & x=y , \end{array} \right.
\ee
\begin{eqnarray*}
& &
\int_{\mathring{G}_1(x)}Z[\vp_2(t)]\tilde{Z}^2(t){\rm d}t-\int_{\mathring{G}_2(y)}Z[\vp_2(t)]\tilde{Z}^2(t){\rm d}t = \\
& &
\frac{2}{\pi}(\sin x+\sin y)H+\mcal{O}(T^{1/6+\epsilon}) .
\end{eqnarray*}
\end{cor}

Next, in the case $x=y=\frac{\pi}{2}$ we obtain

\begin{cor}
\bdis
\begin{split}
& \int_{\mathring{G}_1\left(\frac{\pi}{2}\right)\cup \mathring{G}_2\left(\frac{\pi}{2}\right)}Z[\vp_1(t)]\tilde{Z}^2(t){\rm d}t=\mcal{O}(T^{1/6+\epsilon}), \\
& \int_{\mathring{G}_1\left(\frac{\pi}{2}\right)}Z[\vp_1(t)]\tilde{Z}^2(t){\rm d}t-
\int_{\mathring{G}_2\left(\frac{\pi}{2}\right)}Z[\vp_1(t)]\tilde{Z}^2(t){\rm d}t= \\
& \frac{4}{\pi}H+\mcal{O}(T^{1/6+\epsilon}) ,
\end{split}
\edis
where $[\mathring{T},\widering{T+H}]\subset \mathring{G}_1\left(\frac{\pi}{2}\right)\cup \mathring{G}_2\left(\frac{\pi}{2}\right)$.
\end{cor}

\section{Law of the asymptotic equality of the areas of the positive and negative parts of the graph of the function $Z[\vp_1(t)]\tilde{Z}^2(t)$}

Let
\be \label{3.1}
\begin{split}
\mathring{G}_1^+(x) & =\left\{ t:\ t\in\mathring{G}_1(x),\ Z[\vp_1(t)]\tilde{Z}^2(t)>0\right\} , \\
\mathring{G}_1^-(x) & =\left\{ t:\ t\in\mathring{G}_1(x),\ Z[\vp_1(t)]\tilde{Z}^2(t)<0\right\} , \\
\mathring{G}_2^+(x) & =\left\{ t:\ t\in\mathring{G}_2(x),\ Z[\vp_1(t)]\tilde{Z}^2(t)>0\right\} , \\
\mathring{G}_2^-(x) & =\left\{ t:\ t\in\mathring{G}_2(x),\ Z[\vp_1(t)]\tilde{Z}^2(t)<0\right\} .
\end{split}
\ee
Then we obtain from (\ref{2.1}) by (\ref{1.1}), (\ref{3.1}) the following

\begin{cor}
\be \label{3.2}
\begin{split}
& \int_{\mathring{G}_1^+(x)\cup \mathring{G}_2^+(x)}Z[\vp_1(t)]\tilde{Z}^2(t){\rm d}t \sim \\
& -\int_{\mathring{G}_1^-(x)\cup \mathring{G}_2^-(x)}Z[\vp_1(t)]\tilde{Z}^2(t){\rm d}t, \quad T\to\infty .
\end{split}
\ee
\end{cor}

Indeed, from (\ref{1.5}) by (\ref{3.1}) we have
\bdis
0<(1-\epsilon)\frac{2}{\pi}H\sin x<\int_{\mathring{G}_1(x)}\leq \int_{\mathring{G}_1^+(x)}\leq \int_{\mathring{G}_1^+(x)\cup \mathring{G}_2^+(x)} ,
\edis
and similarly
\bdis
0<(1-\epsilon)\frac{2}{\pi}H\sin x< - \int_{\mathring{G}_1^-(x)\cup \mathring{G}_2^-(x)} .
\edis
Hence, from (\ref{2.1}), $x=y$, we get
\bdis
\int_{\mathring{G}_1^+(x)}+\int_{\mathring{G}_1^-(x)}+\int_{\mathring{G}_2^+(x)}+\int_{\mathring{G}_2^-(x)}=\mcal{O}(T^{1/6+\epsilon})=o(H)
\edis
(see (\ref{1.2}), i.e. (\ref{3.2})).

\begin{remark}
The formula (\ref{3.2}) represents the law of the asymptotic equality of the areas (measures) of the figures which correspond to the positive part and
negative part, respectively, of the graph of the function
\be \label{3.3}
Z[\vp_1(t)]\tilde{Z}^2(t),\ t\in\mathring{G}_1(x)\cup \mathring{G}_2(x)
\ee
with respect to the disconnected sets $\mathring{G}_1^+(x)\cup\mathring{G}_2^+(x),\ \mathring{G}_1^-(x)\cup\mathring{G}_2^-(x)$. This is one of the laws
governing \emph{chaotic} behavior of the positive and negative values of the function (\ref{3.3}).
\end{remark}

\section{Proof of the Theorem}

\subsection{}

Let is remind that
\bdis
\tilde{Z}^2(t)=\frac{{\rm d}\vp_1(t)}{{\rm d}t},\ \vp_1(t)=\frac{1}{2}\vp(t)
\edis
where
\bdis
\tilde{Z}^2(t)=\frac{Z^2(t)}{2\Phi^\prime_\vp[\vp(t)]}=\frac{Z^2(t)}{\left\{ 1+\mcal{O}\left(\frac{\ln \ln t}{\ln t}\right)\right\}\ln t}
\edis
(see \cite{3}, (3.9); \cite{5}, (1.3); \cite{9}, (1.1), (3.1), (3.2)). Thus, the following lemma holds true (see \cite{8}, (2.5); \cite{9}, (3.3)).

\begin{mydef5}
For every integrable function (in the Lebesgue sense) $f(x),\ x\in [\vp_1(T),\vp_1(T+U)]$ the following is true
\be \label{4.1}
\int_T^{T+U}f[\vp_1(t)]\tilde{Z}^2(t){\rm d}t=\int_{\vp_1(T)}^{\vp_1(T+U)}f(x){\rm d}x,\ U\in \left.\left( 0,\frac{T}{\ln T}\right.\right] ,
\ee
where $t-\vp_1(t)\sim (1-c)\pi(t)$.
\end{mydef5}

\begin{remark}
The formula (\ref{4.1}) remains true also in the case when the integral on the right-hand side of eq. (\ref{4.1}) is only relatively convergent improper
integral of the second kind (in the Riemann sense).
\end{remark}

In the case (comp. (\ref{1.7})) $T=\vp_1(\mathring{T}),\ T+U=\vp_1(\widering{T+U})$ we obtain from (\ref{4.1})
\be \label{4.2}
\int_{\mathring{T}}^{\widering{T+U}}f[\vp_1(t)]\tilde{Z}^2(t){\rm d}t=\int_T^{T+U}f(x){\rm d}x .
\ee

\subsection{}

First of all, we have from (\ref{4.2}), for example,
\be \label{4.3}
\int_{\mathring{t}_{2\nu}(-x)}^{\mathring{t}_{2\nu}(x)}f[\vp_1(t)]\tilde{Z}^2(t){\rm d}t=\int_{t_{2\nu}(-x)}^{t_{2\nu}(x)}f(t){\rm d}t ,
\ee
(see (\ref{1.4}). Next, in the case
\bdis
f(t)=Z[\vp_1(t)]\tilde{Z}^2(t),\ t\in \mathring{G}_1(x)\cup \mathring{G}_2(y)
\edis
we have the following $\tilde{Z}^2$-transformation
\be \label{4.4}
\begin{split}
& \int_{\mathring{G}_1(x)}Z[\vp_1(t)]\tilde{Z}^2(t){\rm d}t=\int_{G_1(x)}Z(t){\rm d}t , \\
& \int_{\mathring{G}_2(y)}Z[\vp_1(t)]\tilde{Z}^2(t){\rm d}t=\int_{G_2(y)}Z(t){\rm d}t ,
\end{split}
\ee
(see (\ref{1.1}), (\ref{4.3})). Let us remind that we have proved the following mean-value formulae (see \cite{2})
\be \label{4.5}
\begin{split}
& \int_{G_1(x)}Z(t){\rm d}t=\frac 2\pi H\sin x+\mcal{O}(T^{1/6+\epsilon}) ,  \\
& \int_{G_2(y)}Z(t){\rm d}t=-\frac 2\pi H\sin y+\mcal{O}(T^{1/6+\epsilon}) .
\end{split}
\ee
Now, our formulae (\ref{1.5}) follow from (\ref{4.4}) by (\ref{4.5}).

\thanks{I would like to thank Michal Demetrian for helping me with the electronic version of this work.}

\end{document}